\documentclass[a4 paper,11pt]{amsart}
\usepackage{amsmath,amssymb,amsthm,ascmac}
\usepackage[mathscr,mathcal]{eucal}
\usepackage{latexsym}
\usepackage{enumerate}
\usepackage{graphicx}

\begin{document}

\allowdisplaybreaks

\theoremstyle{definition}
\newtheorem{defi}{\textbf{Definition}}[section]
\newtheorem{thm}[defi]{\textbf{Theorem}}
\newtheorem{lem}[defi]{\textbf{Lemma}}
\newtheorem{prop}[defi]{\textbf{Proposition}}
\newtheorem{cor}[defi]{\textbf{Corollary}}
\newtheorem{ex}[defi]{\textbf{Example}}
\newtheorem{rem}[defi]{\textbf{Remark}}
\newtheorem*{corr}{\textbf{Corollary}}

\theoremstyle{plain}
\newtheorem{maintheorem}{Theorem}
\newtheorem{theorem}{Theorem }[section]
\newtheorem{proposition}[theorem]{Proposition}
\newtheorem{mainproposition}{Proposition}
\newtheorem{lemma}[theorem]{Lemma}
\newtheorem{corollary}[theorem]{Corollary}
\newtheorem{maincorollary}{Corollary}
\newtheorem{claim}{Claim}
\renewcommand{\themaintheorem}{\Alph{maintheorem}}
\theoremstyle{definition} \theoremstyle{remark}
\newtheorem{remark}[theorem]{Remark}
\newtheorem{example}[theorem]{Example}
\newtheorem{definition}[theorem]{Definition}
\newtheorem{problem}{Problem}
\newtheorem{question}{Question}
\newtheorem{exercise}{Exercise}

\newtheorem*{subject}{\UTF{FFFD}\UTF{0693}I}
\newtheorem*{mondai}{Problem}
\newtheorem{lastpf}{\CID{1466}\UTF{FFFD}}
\newtheorem{lastpf1}{proof of theorem2.4}
\renewcommand{\thelastpf}{}

\newtheorem{last}{Theorem}

\renewcommand{\thelast}{}
\renewcommand{\proofname}{\textup{Proof.}}

\renewcommand{\theequation}{\arabic{section}.\arabic{equation}}
\makeatletter
\@addtoreset{equation}{section}

\title[Intrinsic ergodicity for factors of $(-\beta)$-shifts]
{Intrinsic ergodicity for factors of $(-\beta)$-shifts}
\author[M. Shinoda]{Mao Shinoda}
\address{Department of Mathematics, Keio University, Yokohama, 223-8522, Japan}
\email{shinoda-mao@keio.jp}
\author[K. Yamamoto]{Kenichiro Yamamoto}
\address{Department of General Education \\
Nagaoka University of Technology \\
Niigata 940-2188, JAPAN
}
\email{k\_yamamoto@vos.nagaokaut.ac.jp}

\subjclass[2010]{}

\date{}

\maketitle
\large

\begin{abstract}
We show that every subshift factor of a ($-\beta$)-shift is intrinsically ergodic, when $\beta\geq \frac{1+\sqrt{5}}{2}$ and the ($-\beta$)-expansion of $1$ is not periodic with odd period.  Moreover, the unique measure of maximal entropy satisfies a certain Gibbs property. This is an application of the technique established by Climenhaga and Thompson to prove intrinsic ergodicity beyond specification. We also prove that there exists a subshift factor of a ($-\beta$)-shift which is not intrinsically ergodic in the cases other than the above. 
\end{abstract}

\section{Introduction}

The $\beta$-transformations and associated $\beta$-shifts were
introduced by R\'enyi \cite{R} to study expansions with non-integer bases $\beta>1$. More precisely,
a $\beta$-transformation $T_{\beta}\colon [0,1)\to [0,1)$ is defined by
$$T_{\beta}(x):=\beta{x}-\lfloor \beta x\rfloor$$
for every $x\in [0,1)$ where $\lfloor \xi \rfloor$ denotes the largest integer no more than $\xi$
and a $\beta$-shift $\Sigma_{\beta}$ is a subshift consisting of $\beta$-expansions induced by a
$\beta$-transformation.
Dynamical properties of $T_{\beta}$ and $\Sigma_{\beta}$ are very important from viewpoints of
ergodic theory and number theory, so for over half a century,
they have been intensively studied by many authors \cite{CT2,IT,P,PS,R,S,T,W1}.

Recently, to investigate expansions with ``negative" non-integer bases,
Ito and Sadahiro \cite{IS} introduced $(-\beta)$-transformations and associated $(-\beta)$-shifts.
For a real number $\beta>1$, a ($-\beta$)-transformation $T_{-\beta}: (0,1]\rightarrow (0,1]$
is a natural modification of $T_{\beta}$ defined by 
\begin{equation}
\label{trans}
	T_{-\beta}(x)=-\beta x+ \lfloor \beta x \rfloor +1
\end{equation}
for every $x\in (0,1]$. 
Not being the same form of ($-\beta$)-transformations studied by Ito and Sadahiro,
it topologically conjugates with the original one (see \cite[\S1]{LS} for details). 
Similar to the case of $\beta$-shifts,
a ($-\beta$)-shift $\Sigma_{-\beta}$ is a subshift consisting of
($-\beta$)-expansions induced by a ($-\beta$)-transformation.
While the definition of $\Sigma_{-\beta}$ is similar to that of $\Sigma_{\beta}$,
dynamical properties of $\Sigma_{\beta}$ are different from that of $\Sigma_{-\beta}$. For example,
\begin{itemize}
\item
it is clear that the number of fixed points of $\Sigma_{-\beta}$ is strictly greater than that of $\Sigma_{\beta}$
for any $\beta>1$;

\item
it is well-known that $\Sigma_{\beta}$ is topologically mixing for all $\beta>1$, but
we can easily find $\beta>1$ so that $\Sigma_{-\beta}$ is not transitive;

\item
we can find $\beta>1$ so that $\Sigma_{\beta}$ is sofic but $\Sigma_{-\beta}$ is not sofic
(see \cite{LS} for instance).
\end{itemize}
Hence it is worth investigating problems considered for $\beta$-shifts in the case of $(-\beta)$-shifts. 
Recently, dynamical properties of $(-\beta)$-shifts are also studied by many authors from viewpoints of
ergodic theory and number theory \cite{DH,DR, HZY,IS,LS}.

In this paper we investigate uniqueness of measures of maximal entropy for factors of ($-\beta$)-shifts.
Dynamical systems with unique measures of maximal entropy are called {\it intrinsically ergodic}.
It is a fundamental problem to study intrinsic ergodicity
at the interface of ergodic theory and topological dynamics.
Dynamical systems with specification is intrinsically ergodic and it is preserved for their factors \cite{Bow}. 
On the other hand, there are also many dynamical systems without specification,
including $\beta$-shifts and $(-\beta)$-shifts for Lebesgue almost $\beta>1$ \cite{Buz}, but
less is known intrinsic ergodicity for systems beyond specification and their factors.
The present paper is largely motivated by results of  Climenhaga and Thompson \cite{CT}.
In the context of symbolic system, they established a new technique to decide intrinsic ergodicity beyond specification 
(See Section \ref{Preliminaries} for details).
As an application they study $\beta$-shifts and give a positive answer to the open problem posed by
Thomsen in \cite[Problem 28.2]{Boy}: 
every subshift factor of a $\beta$-shift is intrinsically ergodic.
Moreover, the unique measure of maximal entropy is the weak*-limit of 
Dirac measures providing equivalent weights across periodic measures. 
On the other hand, while intrinsic ergodicity for ($-\beta$)-shifts is studied by Liao and Steiner \cite{LS}, 
that for factors of ($-\beta$)-shifts were not known. 
So it is natural to ask whether one can apply Climenhaga and Thompson's technique to ($-\beta$)-shifts. 
Our main results are following.

\begin{maintheorem}
\label{main1}
Let $(\Sigma_{-\beta},\sigma)$ be a $(-\beta)$-shift and
$d_{-\beta}(1)$ be a $(-\beta)$-expansion of $1$ (see \S3 for definitions).
\begin{enumerate}
\item
Suppose that $\beta\ge \frac{1+\sqrt{5}}{2}$ and $d_{-\beta}(1)$ is not periodic with odd period.
Then every subshift factor of $(\Sigma_{-\beta},\sigma)$ is intrinsically ergodic. 
Moreover, the unique measure 
satisfies the Gibbs property on some $\mathcal{G}\subset \mathcal{L}(\Sigma_{-\beta})$
(see Definition \ref{gibbs-def})
and is the weak* limit of $\{\mu_n\}$ with the form
\begin{equation}
\label{periodic}
	\mu_n=\frac{1}{\#Per(n)}\sum_{x\in Per(n)} \delta_x
\end{equation}
where $Per(n)=\{x\in \Sigma_{-\beta}: \sigma^n x=x\}$.

\item
Suppose that either $1<\beta<\frac{1+\sqrt{5}}{2}$ or $d_{-\beta}(1)$ is periodic with odd period.
Then there exists a subshift factor of $(\Sigma_{-\beta},\sigma)$ which is not
intrinsically ergodic.
\end{enumerate}
\end{maintheorem}

\begin{rem}
Theorem \ref{main1} (1) states that not only intrinsic ergodicity for factors of $(-\beta)$-shifts,
but also certain Gibbs property of the unique measure of maximal entropy. These are new results.
For $(-\beta)$-shifts (not factors), the intrinsic ergodicity was already known but even in this case,
our result gives more information on the unique measure:
for the unique measure,
the Gibbs property and representation of the unique measure
as a limit of periodic orbit measure are new results.
\end{rem}

\begin{rem}
\label{rem2}
Theorem \ref{main1} is derived from Theorem \ref{epsilon} stated in the end of \S3.
Theorem \ref{epsilon} states that most of $(-\beta)$-shifts satisfy some weak version of
specification, called ``non-uniform specification" with small obstructions to specification in \cite{C}.
Recently, for a class of subshifts with non-uniform specification, under some additional conditions,
several authors prove various dynamical properties such as
the uniqueness of equilibrium states for every H\"older potentials \cite{CC,CT2},
exponentially mixing property of equilibrium states \cite{C}, and
large deviations \cite{CTY}.
One remaining interesting question is to understand which dynamical properties
hold for $(-\beta)$-shifts and their factors.
\end{rem}

\section{Preliminaries}\label{Preliminaries}
Let $A$ be a finite set and $A^{\mathbb{N}}$ the set of all one-sided infinite sequences on the alphabet $A$
with the standard metric $d(x,y)=2^{-t(x,y)}$, where $t(x,y)=\min\{k\in\mathbb{N}:x_k\not=y_k\}$.
Then the shift map $\sigma\colon A^{\mathbb{N}}\to A^{\mathbb{N}}$ defined by
$\sigma(x)_i=x_{i+1}$ is a continuous map. If a compact subset $X\subset A^{\mathbb{N}}$
satisfies $\sigma(X)\subset X$, then $X$ is called a {\it subshift}. 
For a finite sequence $w\in \bigcup_{n\ge 1}A^n$, denote $w^{\infty}:=ww\cdots\in A^{\mathbb{N}}$
the periodic infinite sequence and by
$[w]:=\{x\in X:x_1\cdots x_{|w|}=w\}$ the \textit{cylinder set} for $w$, where $|w|$ denotes the length of $w$.
The \textit{language} $\mathcal{L}(X)$ is defined by
$$\mathcal{L}(X):=\left\{w\in\bigcup_{n\ge 1}A^n:[w]\not=\emptyset\right\}.$$

For a collection $\mathcal{D}\subset\mathcal{L}(X)$, we define the \textit{entropy} of $\mathcal{D}$ by
$$h(\mathcal{D}):=\limsup_{n\rightarrow\infty}\frac{1}{n}\log\sharp\mathcal{D}_n,$$
where $\mathcal{D}_n:=\{w\in\mathcal{D}:|w|=n\}$, and write $h_{\rm top}(X):=h(\mathcal{L}(X))$.
The \textit{entropy of an invariant measure} $\mu$ is given by
$$h(\mu):=\lim_{n\rightarrow\infty}\frac{1}{n}\sum_{w\in\mathcal{L}_n(X)}-\mu[w]\log \mu[w].$$
Then well-known variational principle \cite[Theorem 8.6]{W2} implies that
$$h_{\rm top}(X)=\sup\{h(\mu):\mu\text{ is invariant}\}.$$
An invariant measure which attains this supremum is called a {\it measure of maximal entropy}.
If such a measure exists and unique, then we say that $(X,\sigma)$ is {\it intrinsically ergodic}.

As we said before, we use Climenhaga and Thompson's technique found in \cite{CT,CT2}
to prove Theorem \ref{main1}.
Before stating their results, we give definitions and notations which are appeared in them.
We start with definitions
 of ``non-uniform" version of
specification and Gibbs property.

\begin{defi}[$(Per)$-Specification]
	Let $\mathcal{L}$ be a language and consider a subset $\mathcal{G}\subset\mathcal{L}$. 
	We say that $\mathcal{G}$ has {\it $(Per)$-specification} with gap size $t\in \mathbb{N}$ if
	\begin{itemize}
		\item 
		for all $m\in \mathbb{N}$ and $w^1, \ldots, w^m\in \mathcal{G}$, 
		there exist $v^1, \ldots, v^{m-1}\in \mathcal{L}$ such that 
		$x:=w^1 v^1 w^2 v^2 \cdots v^{m-1} w^m \in \mathcal{L}$
		and $|v^i|=t$ for all $i$ and
		
		\item
		the cylinder $[x]$ contains a periodic point of period exactly $|x|+t$.
		
	\end{itemize}
\end{defi}

\begin{defi}[Gibbs property]
\label{gibbs-def}
Let $X$ be a subshift and $\mathcal{G}\subset \mathcal{L}(X)$.
An invariant measure $\mu$ satisfies the {\it Gibbs property} on $\mathcal{G}$ if
there exists constants $K>1$ such that
\begin{align*}
	\mu([w])\leq K e^{-nh_{\rm top}(X)}
\end{align*} 
for every $w\in \mathcal{L}(X)_n$ and $n\in \mathbb{N}$ and
\begin{align*}
	\mu([w])\geq K^{-1} e^{-nh_{\rm top}(X)}
\end{align*}
for every $w\in \mathcal{G}_n$ and $n\in \mathbb{N}$. 
\end{defi}

For collections of words $\mathcal{A},\mathcal{B}\subset\mathcal{L}(X)$, denote
$$\mathcal{A}\mathcal{B}:=\{vw\in\mathcal{L}(X):v\in\mathcal{A},w\in\mathcal{B}\}.$$
For a decomposition $\mathcal{L}(X)=\mathcal{A}\mathcal{B}$ of  the language and $M\in\mathbb{N}$, we set
$$\mathcal{A}(M):=\{vw\in\mathcal{L}(X):v\in\mathcal{A},w\in\mathcal{B},|w|\le M\}.$$
Now we state two results found in \cite{CT,CT2}.

\begin{thm}(\cite[Theorem C, Remarks 2.1 and 2.3]{CT2})
\label{CT-1}
Let $(X,\sigma)$ be a subshift whose language $\mathcal{L}(X)$ admits a decomposition
$\mathcal{L}(X)=\mathcal{G}\mathcal{C}^s$, and suppose that the following conditions are satisfied:
\begin{enumerate}[(I)]
\item
$\mathcal{G}(M)$ has $(Per)$-specification for every $M\in \mathbb{N}$;

\item
$h(\mathcal{C}^s)<h_{{\rm top}}(X,\sigma)$.
\end{enumerate}
Then $(X,\sigma)$ is intrinsically ergodic.
Moreover, the unique measure satisfies the Gibbs property on $\mathcal{G}$ and
the sequence of probability measures (\ref{periodic}) converges to the unique measure of maximal entropy.
\end{thm}

\begin{thm}(\cite[Corollary 2.3]{CT},\cite[\S\S3.4]{CT2})
\label{CT-2}
Let $(X,\sigma)$ be a shift space whose language admits a decomposition $\mathcal{G}\mathcal{C}$ satisfying
(I), (II), and let $(\tilde{X},\sigma)$ be a subshift factor of $(X,\sigma)$ such that
$h_{\rm top}(\tilde{X},\sigma)>h(\mathcal{C}^s)$. Then $\tilde{X}$ is intrinsically ergodic.
Moreover, the unique measure satisfies the Gibbs property on some
$\tilde{\mathcal{G}}\subset\mathcal{L}(\tilde{X})$ and
the sequence of probability measures (\ref{periodic})
converges to the unique measure of maximal entropy for $(\tilde{X},\sigma)$.
\end{thm}

\section{$(-\beta)$-shifts}
Let $\beta>1$ and $T_{-\beta}$ be the $(-\beta)$-transformation on $(0,1]$ defined by (\ref{trans}).
For the notational simplicity, we set $A:=\{1,2,\cdots,\lfloor\beta\rfloor +1\}$.
Let $I_1:=(0,\frac{1}{\beta})$, $I_i:=[\frac{i-1}{\beta})$ for $2\le i\le \lfloor \beta\rfloor$, and
$I_{\lfloor \beta \rfloor +1}:=[\frac{\lfloor \beta \rfloor}{\beta},1]$.
For $x\in (0,1]$, we define a one-sided sequence $d_{-\beta}(x)\in  A^{\mathbb{N}}$ by $(d_{-\beta}(x))_i=j$ iff $T_{-\beta}^{i-1}(x)\in I_j$ holds. Then we have
$$x=\frac{-(d_{-\beta}(x))_1}{-\beta}+\frac{-(d_{-\beta}(x))_2}{(-\beta)^2}+\frac{-(d_{-\beta}(x))_3}{(-\beta)^3}+\cdots.$$
We call $d_{-\beta}(x)$ a \textit{$(-\beta)$-expansion} of $x$.
Now, we define
$$\Sigma_{-\beta}:={\rm cl}(\{d_{-\beta}(x)\in A^{\mathbb{N}}:x\in(0,1]\}).$$
Then it is easy to see that $\Sigma_{-\beta}$ is a subshift of $A^{\mathbb{N}}$.
We call $(\Sigma_{-\beta},\sigma)$ a $(-\beta)$-shift.
All $(-\beta)$-shifts are characterized by the alternating order.
The \textit{alternating order} $\prec$ is defined by
$x \prec y$ if there exists $i\in \mathbb{N}$ such that $x_k=y_k$ for all $1\leq k\leq  i-1$ and
$(-1)^i(y_i-x_i)<0$. 
Denote by $x\preceq y$ if $x=y$ or $x\prec y$. 
In \cite{IS}, Ito and Sadahiro proved the following characterization of $(-\beta)$-shifts.

\begin{prop}\cite[Theorem 11]{IS}
\label{alt-beta}
Let $(\Sigma_{-\beta},\sigma)$ be a $(-\beta)$-shift.
\begin{enumerate}
\item
If $d_{-\beta}(1)$ is not periodic with odd period, then
$$\Sigma_{-\beta}=\{x\in A^{\mathbb{N}}:\sigma^k(x)\preceq d_{-\beta}(1),\ n\ge 0\}.$$

\item
If $d_{-\beta}(1)=b_1b_2\cdots$ is periodic with odd period $n$, then
$$\Sigma_{-\beta}=\{x\in A^{\mathbb{N}}:(1b_1\cdots b_{n-1}(b_n-1))^{\infty}\preceq \sigma^k(x)\preceq d_{-\beta}(1),
\ n\ge 0\}.$$
\end{enumerate}
\end{prop}

Now we state Theorem \ref{epsilon}, which plays an important role to prove Theorem \ref{main1}
as we mentioned in Remark \ref{rem2}.

\begin{maintheorem}
\label{epsilon}
Let $(\Sigma_{-\beta},\sigma)$ be a $(-\beta)$-shift.
Suppose that $\beta\ge \frac{1+\sqrt{5}}{2}$ and $d_{-\beta}(1)$ is not periodic with odd period.
Then for any $\epsilon>0$, the language $\mathcal{L}(\Sigma_{-\beta})$ of $\Sigma_{-\beta}$ admits a decomposition
$\mathcal{L}(\Sigma_{-\beta})=\mathcal{G}^{\epsilon}\mathcal{C}^{s,\epsilon}$ satisfying the following properties:
\begin{enumerate}[(I)]
\item
$\mathcal{G^{\epsilon}}(M)$ has $(Per)$-specification for every $M\in\mathbb{N}$;

\item
$h(\mathcal{C}^{s,\epsilon})\le\epsilon$.

\end{enumerate}
\end{maintheorem}

We give a proof of Theorem \ref{epsilon} in \S5.

\section{Graph presentation of $(-\beta)$-shifts}
Let $A=\{1, \ldots, b\}$.
A one-sided sequence ${\bf b}=b_1b_2\cdots\in A^{\mathbb{N}}$ is said to be an {\it alternately shift maximal sequence}
if $b_1=b$ and ${\bf b}$ is greater than or equal to, any of its shifted images in the alternate order, that is,
for each $n\ge 1$, $b_nb_{n+1}\cdots\preceq {\bf b}$.
For an alternately shift maximal sequence ${\bf b}$, we set
$\Sigma_{{\bf b}}:=\{x\in A^{\mathbb{N}}:\sigma^k(x)\preceq {\bf b},k\ge 0\}$. Then we can easily check that
\begin{itemize}
\item
$\Sigma_{{\bf b}}$ is a subshift;

\item
$w\in\mathcal{L}(\Sigma_{{\bf b}})$ if and only if
$w_i\cdots w_{|w|}\preceq b_1\cdots b_{|w|-i+1}$ for each $1\le i\le |w|$.

\item
If $d_{-\beta}(1)$ is not periodic with odd period, then $\Sigma_{d_{-\beta}(1)}=\Sigma_{-\beta}$
(see Proposition \ref{alt-beta} (1) for instance).
\end{itemize}
The aim of this section is to construct a countable labelled directed graph, which presents $\Sigma_{{\bf b}}$
and prove its several properties.

A {\it directed graph} $G$ is a pair $(V,E)$ where $V$ and $E$ are disjoint, finite or countable sets,
together with two maps $i\colon E\to V$ and $t\colon E\to V$. The set $V$ is called the set of vertices and
$E$ is called that of edges. The maps $i$ and $t$ assign to each edge $e\in E$ some pair of vertices $(\alpha,\beta)$
where $e$ starts at vertex $i(e)=\alpha$ and terminates at vertex $t(e)=\beta$.
For $n\ge 1$, we set $E^n(G):=\{e_1\cdots e_n\in E^n:t(e_i)=i(e_{i+1}),1\le i\le n-1\}$ and
$E^{\mathbb{N}}(G)=\{e_1e_2\cdots\in E^{\mathbb{N}}:t(e_i)=i(e_{i+1}),i\ge 1\}$.
For a finite set $A$, a function $\varphi\colon E\to A$ is a labelling of edges of $G$ by symbols from
the alphabet $A$. We call the pair $(G,\varphi)$ a {\it labelled directed graph}.
For $e_1e_2\cdots\in E^{\mathbb{N}}(G)$, we define the labelled walk of $e_1e_2\cdots$
by $\varphi_{\infty}(e_1e_2\cdots)=\varphi(e_1)\varphi(e_2)\cdots\in A^{\mathbb{N}}$.
A finite sequence $w\in A^n$ is
said to be {\it finite labelled path} (in $(G,\varphi)$) if there exists $e_1\cdots e_n\in E^n(G)$ such that
$w=\varphi(e_1)\cdots\varphi (e_n)$. Similarly, we say that a one-sided infinite sequence $x\in A^{\mathbb{N}}$
an {\it infinite labelled path} if there exists $e_1e_2\cdots E^{\mathbb{N}}(G)$ such that $x=\varphi_{\infty}(e_1e_2\cdots)$.
Finally, we set
$A^{\mathbb{N}}(G)={\rm cl}(\varphi_{\infty}(E^{\mathbb{N}}(G))$).
Clearly, $A^{\mathbb{N}}(G)$ is a subshift. We say that $(G,\varphi)$ is a presentation of
a subshift $X$ if $X=A^{\mathbb{N}}(G)$ holds.

In what follows we will construct a countable labelled directed graph $(G_{{\bf b}},\varphi)$,
which is a presentation of $\Sigma_{{\bf b}}$.
For $w\in\mathcal{L}(\Sigma_{{\bf b}})$, we set $k(w):=\max\{k\ge 1:w_{|w|-k+1}\cdots w_{|w|}=b_1\cdots b_k\}$
if such $k\ge 1$ exists and $k(w):=0$ otherwise.
We also define $F^w:=\{x\in\Sigma_{{\bf b}}:wx\in\Sigma_{{\bf b}}\}$ the follower set of
$w\in\mathcal{L}(\Sigma_{{\bf b}})$. It is easy to see that $F^w$ is compact.

\begin{lem}
\label{follower}
Let $w,w\in\mathcal{L}(\Sigma_{{\bf b}})$. If $k(w)=k(w')$, then we have
$F^w=F^{w'}$.
\begin{proof}
We set $k:=k(w)=k(w')$. Without loss of generality, we may assume that $k\ge 1$.
Then we have $w_{|w|-k+1}\cdots w_{|w|}=w'_{|w'|-k+1}\cdots w'_{|w'|}=b_1\cdots b_k$, and
$w'_i\cdots w'_{|w'|}\not=b_1 \cdots b_{|\omega'|-i+1}$ for all $1\leq i\leq |\omega'|-k$.
Let $x\in F^w$. Since $wx\in\Sigma_{{\bf b}}$, we have
$\sigma^{|w'|-k}(w'x)=b_1\cdots b_kx=\sigma^{|w|-k}(wx)\in\Sigma_{{\bf b}}$.
This implies that $\sigma^n(w'x)\preceq {\bf b}$ for $n\ge |w'|-k$.
Since $w_i'\cdots w'_{|w'|}\not =b_1\cdots b_{|\omega'|-i+1}$ for all $1\leq i \leq |\omega'|-k$, we have
$\sigma^n(w'x)=w'_{n+1}\cdots w'_{|w'|}x \prec {\bf b}$ for $0\le n\le |w'|-k-1$.
Thus we have $w'x\in\Sigma_{{\bf b}}$, which implies that $F^w\subset F^{w'}$.
$F^{w'}\subset F^w$ can be shown similarly, and so  we have $F^w=F^{w'}$.
\end{proof}
\end{lem}

Now, we define a countable labelled directed graph $(G_{{\bf b}},\varphi)$ as follows.
Let $V=\{V_0,V_1,V_2, \ldots\}$ be a set of vertices.
There is an edge from $V_i$ to $V_j$ whenever there exist $w\in\mathcal{L}(\Sigma_{{\bf b}})$ and $a\in A$
such that $wa\in\mathcal{L}(\Sigma_{{\bf b}})$, $k(w)=i$ and $k(wa)=j$, and this edge is labelled with the symbol $a$.
The following proposition is clear by the definition of $(G_{{\bf b}},\varphi)$ and so we omit the proof.

\begin{prop}
\label{graph1}
Let $(G_{{\bf b}},\varphi)$ be as above.\\
(1) Let $w\in\mathcal{L}(\Sigma_{{\bf b}})$ with $k(w)=k$. Then we have
$$F^w=\{x\in A^{\mathbb{N}}:\text{$x$ is an infinite labelled path started at $V_k$}\}.$$
(2) We have $\Sigma_{{\bf b}}=\{x\in A^{\mathbb{N}}:\text{$x$ is an infinite labelled path}\}$\\
\hspace{78pt}$=\{x\in A^{\mathbb{N}}:\text{$x$ is an infinite labelled path started at $V_0$}\}.$\\
(3) We have
$\mathcal{L}(\Sigma_{{\bf b}})=\{w:\text{$w$ is a finite labelled path}\}$\\
\hspace{95pt}$=\{w:\text{$w$ is a finite labelled path started at $V_0$}\}.$\\
(4) There is an edge from $V_i$ to $V_j$ which is labelled with the symbol $a\in A$ if and only if either
\begin{itemize}
\item
$j=i+1$ and $a=b_{i+1}$ $(i\ge 0)$;

\item
$i$ is odd, $b_{i+1}+1\le a\le b$, $b_1\cdots b_ia\in\mathcal{L}(\Sigma_{{\bf b}})$
and $k(b_1\cdots b_i  a)=j$;

\item
$i$ is even, $1\le a\le b_{i+1}-1$, $b_1\cdots b_ia\in\mathcal{L}(\Sigma_{{\bf b}})$ and
$k(b_1\cdots b_ka)=j$.
\end{itemize}
\end{prop}

\begin{figure}[htbp]
 \begin{center}
  \includegraphics[width=100mm, angle=90]{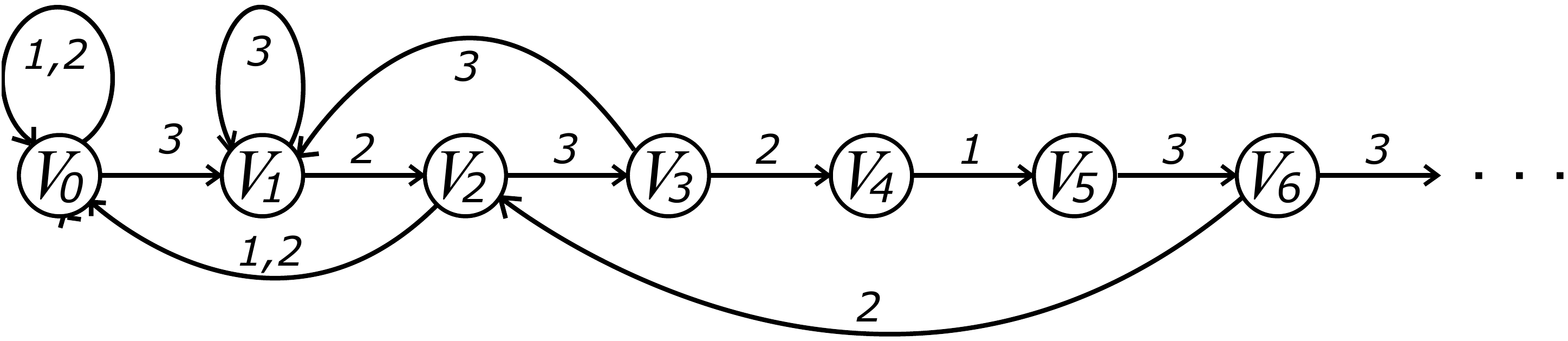}
 \end{center}
 \caption{A graph $(G_{\bf b},\varphi)$ for ${\bf b}=3232133\cdots$.}
 \label{fig-1}
\end{figure}

The situation is sketched in Figure \ref{fig-1} for ${\bf b}=3232133\cdots$.
The following proposition plays a key role to prove Theorem \ref{epsilon} in \S5.

\begin{prop}
\label{entropy-small}
Let $(G_{{\bf b}},\varphi)$ be as above.
For $L\ge 1$, we set
$$\mathcal{C}^{(L)}:=\left\{b_Lw \middle|
\begin{array}{l}
w\text{ is a finite path started at }V_L \\
\text{and never visits }V_0,V_1,\ldots,V_{L-1}
\end{array}
\right\}.$$
Then for any $\epsilon>0$, there exists an $L\ge 1$ such that
$h(\mathcal{C}^{(L)})\le \epsilon$ holds.
\begin{proof}
Let $\epsilon>0$ and take an integer $N\ge 1$ so large that
$\frac{2}{N}\log b+\frac{2}{N}\log N\le \epsilon$ holds.
First, we prove the following lemma.
\begin{lem}
\label{keyprop-1}
There exists $L\ge 1$ such that
for any $k\ge L$,
$$\{e\in E:i(e)=V_k,\ t(e)=V_j,\ 0\le k-j\le N\}=\emptyset.$$
\begin{proof}
By contradiction, assume that
for any $L\ge 1$, we can find $k\ge L$ such that
$$\{e\in E:i(e)=V_k,\ t(e)=V_j,\ 0\le k-j\le N\}\not=\emptyset.$$
Then we can find a strictly increasing sequence $\{k_i\}_{i\ge 1}$ such that for any $i\ge 1$, $k_{i+1}-k_i>N$ and
$$\{e\in E:i(e)=V_{k_i},\ t(e)=V_j,\ 0\le k_i-j\le N\}\not=\emptyset$$
hold. So there exist $\{e_i\}_{i\ge 1}\subset E$ and $\{j_i\}_{i\ge 1}\subset\mathbb{N}$ such that for any $i\ge 1$,
$i(e_i)=V_{k_i}$, $t(e_i)=V_{j_i}$ and $0\le k_i-j_i\le N$ hold.
Thus, we can find $i$ and $i'$ such that $i>i'$ and $k_i-j_i=k_{i'}-j_{i'}$.
Note that $j_i\ge k_i-N>k_{i'}\ge j_{i'}$.
Since $i(e_p)=V_{k_p}$ and $t(e_p)=V_{j_p}$ for $p=i,i'$,
by the definition of the graph $G_{{\bf b}}$, we have
$$b_{k_i-j_i+2}\cdots b_{k_i}\varphi(e_i)=b_1\cdots b_{j_i}\text{ and}$$
$$b_{k_{i'}-j_{i'}+2}\cdots b_{k_{i'}}\varphi(e_{i'})=b_1\cdots b_{j_{i'}}.$$
Therefore, we have
$$b_{k_{i'}+1}\not =\varphi(e_{i'})=b_{j_{i'}}=b_{j_{i'}+k_i-j_i+1}=b_{k_{i'}+1},$$
which is a contradiction.
\end{proof}
\end{lem}
It follows from Lemma \ref{keyprop-1} that we can find $L\ge 1$ such that
for any $k\ge L$,
\begin{equation}
\label{far-edge}
\{e\in E:i(e)=V_k,\ t(e)=V_j,\ 0\le k-j\le N\}=\emptyset
\end{equation}
holds. In what follows we will show that $h(\mathcal{C}^{(L)})\le \epsilon$.
We set
$$V':=\{V_{L-1},V_L,V_{L+1},\ldots\},\ E':=\{e\in E:i(e)\in V',t(e)\in V'\setminus\{V_{L-1}\}\}$$
and define a subgraph $G'=(V',E')$ of $G_{{\bf b}}$.
For the notational simplicity, we set $W_i:=V_{L+i-2}$ $(i\ge 1)$.
Let $A=(a_{ij})$ be a adjacency matrix of $G'$, that is,
$a_{ij}=\sharp\{e\in E':i(e)=W_i,t(e)=W_j\}$.
We denote by $A^n=(a_{ij}^{(n)})$ for $n\ge 1$ and set
$a_i^{(n)}:=\sum_{j=1}^{\infty}a_{ij}^{(n)}$. Then it is easy to see that
$a_1^{(n)}=\sharp\mathcal{C}^{(L)}_n.$
Thus, to prove the proposition, it is sufficient to show that
$$\limsup_{n\rightarrow\infty}\frac{1}{n}\log a_1^{(n)}\le\epsilon.$$
By the equation (\ref{far-edge}) and the definition of the graph $G_{{\bf b}}$,
the adjacency matrix $A$ and $a_i^{(n)}$ satisfies the following properties.
\begin{itemize}
\item
$a_{ij}\in\{0,1\}$ for any $i,j\ge 1$.

\item
$a_{i,i+1}=1$ for any $i\ge 1$.

\item
$j\not=i+1$ and $i-N< j$, then $a_{ij}=0$.

\item
$a_i^{(1)}=1$ for $1\le i\le N$ and $a_i^{(1)}\le b$ for any $i\ge 1$.

\item
$a_i^{(n)}\le a_i^{(n+1)}$ for any $n\ge 1$ and $i\ge 1$.
\end{itemize}
Thus, for any $n\ge 1$, we have the following estimates of $a_i^{(n+1)}$.
\begin{itemize}
\item
If $1\le i\le N$, then $a_i^{(n+1)}\le a_{i+1}^{(n)}$.

\item
If $i\ge N+1$, then we can find $(s_1,\ldots,s_{i-N})\in\{0,1\}^{i-N}$ such that
$\sharp\{1\le j\le i-N:s_j=1\}\le b-1$ and 
$$a_i^{(n+1)}\le \left(\sum_{j=1}^{i-N}s_ja_j^{(n)}\right)+a_{i+1}^{(n)}.$$
\end{itemize}
In what follows we will estimate $a_1^{(qN+1)}$ for $q\ge 2$.
\vspace{0.1cm}\\
({\em Step} 1) $a_1^{(2N+1)}\le bN$.\\
\begin{eqnarray*}
{}
& &
a_1^{(2N+1)}\\
&\le&
a_2^{(2N)}\le\cdots\le a_{N+1}^{(N+1)}\\
&\le&
a_1^{(N)}+a_{N+2}^{(N)}\\
&\le&
1+a_{N+2}^{(N)}\\
&\le&
1+(a_1^{(N-1)}+a_2^{(N-1)})+a_{N+3}^{(N-1)}\\
&\le&
(1+2)+a_{N+3}^{(N-1)}\\
&\le&
\cdots\\
&\le&
\frac{(b-1)b}{2}+a_{N+b}^{(N-b+2)}\\
&\le&
\frac{(b-1)b}{2}+\left(\sum_{j=1}^bs_ja_j^{(N-b+1)}\right)+a_{N+b+1}^{(N-b+1)}\\
& & (\text{for some }(s_1,\ldots, s_b)\in\{0,1\}^b\text{ with }\sharp\{1\le j\le b:s_j=1\}\le b-1)\\
&\le&
\frac{(b-1)b}{2}+(b-1)+a_{N+b+1}^{(N-b+1)}\\
&\le&
\frac{(b-1)b}{2}+2(b-1)+a_{N+b+2}^{(N-b)}\\
&\le&
\cdots\\
&\le&
\frac{(b-1)b}{2}+(N-b+1)(b-1)+a_{2N+1}^{(1)}\\
&\le&
\frac{(b-1)b}{2}+(N-b+1)(b-1)+b\\
&\le&
bN.
\end{eqnarray*}
By the above estimation, we also have
$a_i^{(2N+2-i)}\le bN$ for any $N+1\le i\le 2N+1$.
\vspace{0.1cm}\\
({\em Step} 2) $a_1^{(3N+1)}\le b^3N^3$.\\
\begin{eqnarray*}
{}
& &
a_1^{(3N+1)}\\
&\le&
a_2^{(3N)}\le\cdots\le a_{N+1}^{(2N+1)}\\
&\le&
a_1^{(2N)}+a_{N+2}^{(2N)}\\
&\le&
bN+a_{N+2}^{(N)}\\
&\le&
bN+(a_1^{(2N-1)}+a_2^{(2N-1)})+a_{N+3}^{(2N-1)}\\
&\le&
(1+2)bN+a_{N+3}^{(2N-1)}\\
&\le&
\cdots\\
&\le&
\frac{(b-1)b^2N}{2}+a_{N+b}^{(2N-b+2)}\\
&\le&
\frac{(b-1)b^2N}{2}+\left(\sum_{j=1}^bs_ja_j^{(2N-b+1)}\right)+a_{N+b+1}^{(2N-b+1)}\\
& & (\text{for some }(s_1,\ldots, s_b)\in\{0,1\}^b\text{ with }\sharp\{1\le j\le b:s_j=1\}\le b-1)\\
&\le&
\frac{(b-1)b^2N}{2}+(b-1)bN+a_{N+b+1}^{(2N-b+1)}\\
&\le&
\frac{(b-1)b^2N}{2}+2(b-1)bN+a_{N+b+2}^{(2N-b)}\\
&\le&
\cdots\\
&\le&
\frac{(b-1)b^2N}{2}+(N-b+1)(b-1)bN+a_{2N+1}^{(N+1)}.
\end{eqnarray*}
Then,
\begin{eqnarray*}
{}
& &
a_{2N+1}^{(N+1)}\\
&\le&
\left(\sum_{j=1}^{N+1}t_ja_j^{(N)}\right)+a_{2N+2}^{(N)}\\
& & (\text{for some }(t_1,\ldots, t_{N+1})\in\{0,1\}^b\text{ with }\sharp\{1\le j\le b:t_j=1\}\le b-1)\\
&\le&
(b-1)bN+a_{2N+2}^{(N)}\\
&\le&
2(b-1)bN+a_{2N+3}^{(N-1)}\\
&\le&
\cdots\\
&\le&
(b-1)bN^2+a_{3N+1}^{(1)}\\
&\le&
(b-1)bN^2+b\\
&\le&
b^2N^2.
\end{eqnarray*}
Therefore, we have
$$a_1^{(3N+1)}\le\frac{(b-1)b^2N}{2}+(N-b+1)(b-1)bN+b^2N^2\le b^3N^3.$$
By above estimations, we also have
$$
a_i^{(3N+2-i)}\le
\left\{
\begin{array}{ll}
b^3N^3 & (N+1\le i\le 2N); \\
b^2N^2 & (2N+1\le i\le 3N).
\end{array}
\right.
$$
Using this procedure inductively, we have
$$a_1^{(qN+1)}\le b^{2q-3}N^{2q-3}\text{ for any }q\ge 2.$$
Fix any $n\ge N$. Then we can find $q\ge 1$ and $0\le r\le N-1$ such that
$n=qN+r$. Since $n\le (q+1)N+1$, we have
$$a_1^{(n)}\le a_1^{(q+1)N+1}\le b^{2q-1}N^{2q-1}.$$
Therefore, we have
\begin{eqnarray*}
\frac{1}{n}\log a_1^{(n)}
&\le&
\frac{1}{n}\log b^{2q-1}+\frac{1}{n}\log N^{2q-1}\\
&\le&
\frac{2}{N}\log b+\frac{2}{N}\log N \\
&\le&
\epsilon,
\end{eqnarray*}
which proves the proposition.
\end{proof}
\end{prop}

\section{Proof of Theorem \ref{epsilon}}

Throughout of this section, let $\beta\ge \frac{1+\sqrt{5}}{2}$ and $d_{-\beta}(1)$ is
not periodic with odd period.
Let $T_{-\beta}\colon (0,1]\to (0,1]$ be a $(-\beta)$-transformation.
We extend $T_{-\beta}$ to $[0,1]$ by
$$
T_{-\beta}(x)=
\left\{
\begin{array}{ll}
T_{-\beta}(x) & (x\in (0,1]); \\
1 & (x=0).
\end{array}
\right.
$$
Then it is known that $T_{-\beta}$ is locally eventually onto in the following sense.

\begin{lem}(\cite[Theorem 2.2]{LS})
\label{onto}
For any interval $I\subset [0,1]$ of positive length, we can find an integer $n>0$ such that
$T_{-\beta}^n(I)=(0,1]$.
\end{lem}

We define a map $\Psi\colon \Sigma_{-\beta}\to [0,1]$ by
$\{\Psi(x)\}=\bigcap_{n=0}^{\infty}T_{-\beta}^{-n}({\rm cl}(I_{x_{n+1}}))$.
Since $d_{-\beta}(1)$ is not periodic with odd period, $\Psi$ satisfies the following properties.
\begin{itemize}
\item
$\Psi$ is continuous, surjective, and $\Psi\circ\sigma=T_{-\beta}\circ\Psi$.

\item
 $\Psi^{-1}\{x\}$ is a singleton for every $x\in [0,1]\setminus\bigcup_{n=0}^{\infty}T_{-\beta}^{-n}(\{1\})$

\item
For any $w\in\mathcal{L}(\Sigma_{-\beta})$, $\Psi([w])$ is a closed interval with positive length.
\end{itemize}

\begin{lem}
\label{mixing}
$(\Sigma_{-\beta},\sigma)$ is topologically mixing.
\begin{proof}
Let $v,w\in\mathcal{L}(\Sigma_{-\beta})$.
Then $\Psi([v])$ and $\Psi([w])$ is an interval of positive length.
It follows from Lemma \ref{onto} that there is $N>0$ such that
$\Psi(\sigma^n([v]))=T_{-\beta}^n(\Psi([v]))=(0,1]$ for any $n\ge N$, which implies that
$\Psi(\sigma^n([v]))\cap\Psi([w])$ is an interval of positive length.
Note that the inverse image of a point by $\Psi$ is a singleton
except at most countable set.
So we have $\sigma^n([v])\cap [w]\not=\emptyset$, which proves the lemma.
\end{proof}
\end{lem}

In the rest of this section, denote ${\bf b}=d_{-\beta}(1)$, and let $(G_{{\bf b}},\varphi)$
be a labelled directed graph as in \S4. We denote $V=\{V_0,V_1,\ldots\}$ the set of vertices and
$E$ the set of edges.

\begin{lem}
\label{irreducible}
Suppose that $d_{-\beta}(1)$ is not eventually periodic.
Then for any $i\ge 1$, we can find a finite path from $V_i$ to $V_0$.
\begin{proof}
Since $\Psi([b_1\cdots b_i])$ is an interval with positive length, it follows from Lemma \ref{onto}
that there is an integer $n>0$ such that $\Psi(\sigma^n([b_1\cdots b_i]))=(0,1]$.
Thus, we have $\Psi(1^{\infty})\in\Psi(\sigma^n([b_1\cdots b_i]))$.
Since $\Psi$ is injective at $1^{\infty}$, we have $1^{\infty}\in\sigma^n([b_1\cdots b_i])$.
So we can find $v\in\mathcal{L}_{n-i}(\Sigma_{-\beta})$ such that $b_1\cdots b_i v1^{\infty}\in\Sigma_{-\beta}$.
Since $k(b_1\cdots b_i)=i$, by the definition of $G_{{\bf b}}$, there exists an infinite path
$e_1e_2\cdots\in E^{\mathbb{N}}$ such that $i(e_1)=V_i$ and $\varphi_{\infty}(e_1e_2\cdots)=v1^{\infty}$.

Thus, it is sufficient to show that there exists an integer $j\ge 1$ such that $t(e_j)=V_0$.
By contradiction, assume that $t(e_j)\not=V_0$ for any $j\ge 1$. Since $d_{-\beta}(1)$ is not eventually periodic,
we can find $k\ge i$ such that $b_k\not=1$. Note that $\varphi(e_j)=1$ for any $j\ge n+1$.
Thus for each $j\ge 1$,
we can find an integer $1\le i(j)\le k$ so that $t(e_j)=V_{i(j)}$ holds.
Hence we have $1\le k(v1^{2k})\le k$, which implies that $b_1=1$.
This is a contradiction.
\end{proof}
\end{lem}

In what follows we will prove Theorem \ref{epsilon}.\\
({\em Case} 1) $d_{-\beta}(1)$ is eventually periodic.\\
In this case, it follows from \cite[Theorem 12]{IS} that $(\Sigma_{-\beta},\sigma)$ is sofic.
Thus, it follows from Lemma \ref{mixing} that $\mathcal{L}(\Sigma_{-\beta})$ has (Per)-specification.
Hence the conclusion of Theorem \ref{epsilon} is trivial.
\vspace{0.1cm}\\
({\em Case} 2) $d_{-\beta}(1)$ is not eventually periodic.\\
Take any $\epsilon>0$. By Proposition \ref{entropy-small}, we can find $L>0$ such that
$h(\mathcal{C}^{(L)})\le \epsilon$. Here we set
$$\mathcal{C}^{(L)}:=\left\{b_Lw \middle|
\begin{array}{l}
w\text{ is a finite path started at }V_L \\
\text{and never visits }V_0,V_1,\ldots,V_{L-1}
\end{array}
\right\}.$$
Now, we define
$$\mathcal{G}^{(L)}:=\left\{w \middle|
\begin{array}{l}
w\text{ is a finite path started at }V_0 \\
\text{and ended at }V_i\text{ for some }0\le i\le L-1
\end{array}
\right\}$$
and set $\mathcal{G}^{\epsilon}:=\mathcal{G}^{(L)}$ and $\mathcal{C}^{s,\epsilon}:=\mathcal{C}^{(L)}$.
Clearly, we have $\mathcal{L}(\Sigma_{-\beta})=\mathcal{G}^{\epsilon}\mathcal{C}^{s,\epsilon}$
and $h(\mathcal{C}^{s,\epsilon})\le\epsilon$.
Finally, we prove that $\mathcal{G}^{\epsilon}(M)$ has (Per)-specification for every $M\in\mathbb{N}$.
By Lemma \ref{irreducible}, we can define
$$t_M:=\max\{\text{length of shortest path from $V_i$ to }V_0:0\le i\le M+L-1\}.$$
Let $m\in\mathbb{N}$ and $w^1,\ldots,w^m\in\mathcal{G}^{\epsilon}(M)$.
Then for each $1\le i\le n$, $w^i$ is a finite labelled path started at $V_0$ and ended at $V_{k_i}$ for
some $0\le k_i\le M+L-1$.
Thus by the definition of  $t_M$,
we can find a finite labelled path $u^i$ started at $V_{k_i}$ and ended at $V_0$ so that $|u^i|\le t_M$ holds.
If we set $v^i:=u^i1^{t_M-|u^i|}$, then it is easy to see that
$|v^i|=t_M$ and $w^iv^i$ is a finite labelled path started at $V_0$ and ended at $V_0$.
Thus we have $(w^1v^1w^2v^2\cdots w^mv^m)^{\infty}\in \Sigma_{-\beta}$, which implies that
$\mathcal{G}^{\epsilon}(M)$ has (Per)-specification.
Theorem \ref{epsilon} is proved.

\begin{rem}
It is known that $\beta$-shifts $\Sigma_{\beta}$ also can be presented by a countable labelled directed graph.
In this case, a decomposition of language $\mathcal{L}(\Sigma_{\beta})=\mathcal{G}^{(1)}\mathcal{C}^{(1)}$ satisfies
$h(\mathcal{C}^{(1)})=0$ \cite[\S\S3.1]{CT} and so Theorem \ref{CT-2} can be applicable to
factors of $\beta$-shifts.
In contrast, in the case of $(-\beta)$-shifts, the entropy of $\mathcal{C}^{(1)}$ is large for most of $\beta$s,
and it seems to be hard to find $L$ so that $h(\mathcal{C}^{(L)})=0$, which makes the
proof of Theorem \ref{main1} difficult. This difficulty comes from the complexity of the graph presentation of
$(-\beta)$-shifts (see Figure \ref{fig-1}).
The novelty of this paper is to overcome this difficulty by showing $\lim_{L\rightarrow\infty}h(\mathcal{C}^{(L)})=0$
(Proposition \ref{entropy-small}),
which allows us to use Theorem \ref{CT-2} to factors of $(-\beta)$-shifts.
\end{rem}

\section{Proof of Theorem \ref{main1}}

\noindent
(1) Suppose that $\beta\ge\frac{1+\sqrt{5}}{2}$ and $d_{-\beta}(1)$ is not periodic with odd period.
Let $(X,\sigma)$ be a subshift factor of $(\Sigma_{-\beta},\sigma)$.\\
({\em Case} 1) $h_{{\rm top}}(X,\sigma)=0$.\\
It follows from \cite[Propositions 2.2 and 2.4]{CT} that $X$ comprises a single periodic orbit
and so $(X,\sigma)$ is intrinsically ergodic.
\vspace{0.1cm}\\
({\em Case} 2) $h_{{\rm top}}(X,\sigma)>0$.\\
Take $\epsilon>0$ so small that $\epsilon<h_{{\rm top}}(X,\sigma)$ holds.
Then it follows from Theorem \ref{epsilon} that $\mathcal{L}(\Sigma_{-\beta})$ admits a
decomposition $\mathcal{L}(\Sigma_{-\beta})=\mathcal{G}^{\epsilon}\mathcal{C}^{s,\epsilon}$ satisfying
\begin{enumerate}[(I)]
\item
$\mathcal{G^{\epsilon}}(M)$ has $(Per)$-specification for every $M\in\mathbb{N}$;

\item
$h(\mathcal{C}^{s,\epsilon})\le\epsilon<h_{{\rm top}}(X,\sigma)$;

\end{enumerate}

Thus it follows from Theorem \ref{CT-2} that $(X,\sigma)$ is intrinsically ergodic.
\vspace{0.3cm}\\
(2) Suppose that either $1<\beta<\frac{1+\sqrt{5}}{2}$ or $d_{-\beta}(1)$ is periodic with odd period.
Let $X:=\{1^{\infty}\}\cup\{1^n2^{\infty}:n\ge 1\}\cup\{2^{\infty}\}$.
Then it is easy to see that $(X,\sigma)$ is a subshift with zero topological entropy, and
has exactly two ergodic measures of maximal entropy $\delta_{1^{\infty}}$, $\delta_{2^{\infty}}$.
In particular, $(X,\sigma)$ is not intrinsically ergodic.
So it is sufficient to show that there is a continuous surjection $\pi\colon\Sigma_{-\beta}\to X$
such that $\pi\circ\sigma=\sigma\circ\pi$ holds.\\
({\em Case} 1) $1<\beta<\frac{1+\sqrt{5}}{2}$.\\
Since $\beta<\frac{1+\sqrt{5}}{2}$, we have
$d_{-\beta}(1)<d_{-\beta'}(1)=21^{\infty}$. Here we set $\beta':=\frac{1+\sqrt{5}}{2}$.
So $d_{-\beta}(1)$ is of the form ``$21^k2\cdots$" for some even number $k$.
We note that $k\ge 2$ (see \cite[Proof of Proposition 3.5]{LS} for instance).
Let $n:=k+1$. Recall that $w\in\mathcal{L}(\Sigma_{-\beta})$ if and only if
$$w_i\cdots w_{|w|}\le b_1\cdots b_{|w|-i+1}\ (1\le i\le |w|).$$
Thus, we have
\begin{equation}
\label{not-in}
w1^n\not\in\mathcal{L}(\Sigma_{-\beta})\text{ whenever }w\not= 1^{|w|}.
\end{equation}
We define a block map $\Phi\colon\mathcal{L}_n(\Sigma_{-\beta})\to\{1,2\}$ as
$$\Phi(w)=
\left\{
\begin{array}{ll}
1 & (w=1^n); \\
2 & (\text{otherwise}),
\end{array}
\right.
$$
and define $\pi\colon\Sigma_{-\beta}\to X$ by
$$(\pi(x))_k:=\Phi(x_k\cdots x_{k+n-1})\ (k\ge 1).$$
First, we show that $\pi$ is well-defined, that is, $\pi(x)\in X$ for any $x\in\Sigma_{-\beta}$.
%
If $\pi(x)=1^{\infty}$, then clearly $\pi(x)\in X$.
Assume $\pi(x)\not=1^{\infty}$ and let $k\ge 1$ be a minimum integer so that $(\pi(x))_k=2$.
Then it is sufficient to show that $(\pi(x))_j=2$ for any $j\ge k+1$. 
Since $(\pi(x))_k=2$, by the definition of a block map, we have
$x_k\cdots x_{n-k+1}\not=1^n$.
Now, by contradiction, assume that $(\pi(x))_j=1$ for some $j\geq k+1$. 
Then $x_j\cdots x_{n+j-1}=1^n$. This together with $x_k\cdots x_{n-k+1}\not=1^n$ implies that 
$x_k\cdots x_{j-1} \not=1^{j-k}$. Thus, by the equation (\ref{not-in}), we have
$x_k\cdots x_{j-1}x_j\cdots x_{n+j-1} \notin \mathcal{L}(\Sigma_{-\beta})$, which is contradiction. 
Therefore $\pi\colon\Sigma_{-\beta}\to X$ is well-defined.

Clearly $\pi$ is a continuous and satisfies $\pi\circ\sigma=\sigma\circ\pi$.
Finally, we show that $\pi$ is surjective.
As we proved above, $(\pi(x))_k=2$ implies that $(\pi(x))_j=2$ for any $j\ge k$.
Since $\pi(d_{-\beta}(1))_1=\pi(1^{n-1}d_{-\beta}(1))=2$, we have $\pi(d_{-\beta}(1))=2^{\infty}$ and
$\pi(1^{k+n}d_{-\beta}(1))=1^k2^{\infty}$ for any $k\ge 1$.
These together with $\pi(1^{\infty})=1^{\infty}$ imply that $\pi$ is surjective.
\vspace{0.3cm}\\
({\em Case} 2) $d_{-\beta}(1)$ is periodic with odd period $n$.\\
First, we recall that 
\begin{align*}
	\Sigma_{-\beta}=\{x: (1b_1\cdots b_{n-1}(b_n-1))^{\infty}\preceq \sigma^k(x) \preceq d_{-\beta}(1) ,k\ge 1 \},
\end{align*}
(see Proposition \ref{alt-beta}).
We claim that $[b_i\cdots b_{i+3n-1}]=\{(b_i\cdots b_{i+n-1})^{\infty}\}$ for $1\le i\le n$.
By contradiction, assume that there exists $x\in [b_i\cdots b_{i+3n-1}]$ so that $x\not =(b_i\cdots b_{i+n-1})^{\infty}$.
Then we have $y:=\sigma^{n-i+1}(x)\in [b_1\cdots b_{2n}]$ and $y\not=(b_1\cdots b_n)^{\infty}$.
Let $k>2n$ be a smallest number so that $y_k\not=b_k$ holds. It follows from $y\prec d_{-\beta}(1)$ that
$(-1)^k(b_k-y_k)<0$. Since $n$ is a odd number, we have $(-1)^{k-n}(b_{k-n}-y_k)>0$, which implies that
$\sigma^n(y)\succ d_{-\beta}(1)$. This contradicts with $\sigma^n(y)\in\Sigma_{-\beta}$.
We define a block map $\Phi\colon\mathcal{L}_{3n}(\Sigma_{-\beta})\to\{1,2\}$ as
$$\Phi(w)=
\left\{
\begin{array}{ll}
2 & (w=b_i\cdots b_{i+3n-1},\ 1\le i\le n); \\
1 & (\text{otherwise}),
\end{array}
\right.
$$
and define $\pi\colon\Sigma_{-\beta}\to X$ by
$$(\pi(x))_k:=\Phi(x_k\cdots x_{k+3n-1})\ (k\ge 1).$$

In what follows we will show that $\pi$ is a factor map.
First, we show that $(\pi(x))_k=2$ for some $k\ge 1$ implies that $(\pi(x))_j=2$ for $j\geq k$.
Assume $(\pi(x))_k=2$. 
Then there exists $1\leq i\leq n$ such that
\begin{align*}
	x_k\cdots x_{k+3n-1}=b_i\cdots b_{i+3n-1}.
\end{align*}
Hence $x_k x_{k+1}\cdots \in [b_i\cdots b_{i+3n-1}]=\{(b_i\cdots b_{i+n-1})^\infty\}$
and we have $(\pi(x))_j=2$.
Therefore, $\pi$ is well-defined.
It is clear that $\pi$ is continuous and $\pi\circ\sigma=\sigma\circ\pi$.

Finally, we prove the surjectivity of $\pi$. 
It is easy to see that $\pi(1^\infty)=1^\infty$ and 
$\pi(d_{-\beta}(1))=2^\infty$.
Fix $k\geq 1$.
It is enough to show that there exists $x\in \Sigma_{-\beta}$ such that $\pi(x)=1^k 2^\infty$. 
Note that $b_n\neq 1$. Indeed, since $T_{-\beta}^n(1)=1$, we have $T_{-\beta}^{n-1}(1)=\frac{p}{\beta}$
for some $1\le p\le \lfloor\beta\rfloor$. This implies that $b_n\not=1$ by the definition of $\beta$-expansions.
We claim that for all $j\geq1$ and $1\leq i \leq n$ 
\begin{equation}
\label{final}
	1^j b_1\cdots b_{3n-j} \neq b_i\cdots b_{i+3n-1}.
\end{equation}
By contradiction, assume that $1^jb_1\cdots b_{3n-j}=b_i\cdots b_{i+3n-1}$.
Then we have $1^jd_{-\beta}\in[b_i\cdots b_{i+3n-1}]=(b_i\cdots b_{i+n-1})^{\infty}$.
Hence we have $b_n=1$, which is a contradiction.
By the equation (\ref{final}), we have $\pi(1^kd_{-\beta}(1))=1^k2^{\infty}$.
Therefore, $\pi$ is surjective.

\vspace*{3mm}

\noindent
\textbf{Acknowledgement.}~ 
The first author was partially supported by Grant-in-Aid for JSPS Research Fellow of the JSPS, 17J03495.
The second author was partially supported by JSPS KAKENHI Grant Number 18K03359.

\end{document}